\documentclass{article}
\newcommand{\dss}{\displaystyle}

\usepackage[psamsfonts]{amssymb}
\usepackage{amsmath,amsthm}
\usepackage{amsopn}

\begin{document}

\setlength{\baselineskip}{1.75em}

\title{Change of Variable for Multi-dimensional Integral}
\date{4 March 2003}

\maketitle

\centerline{Isidore Fleischer}

\vskip.110in

\abstract{The change of variable theorem is proved under the sole
hypothesis of differentiability of the transformation.
Specifically, it is shown under this hypothesis that the
transformed integral equals the given one over every
measurable subset on which the transformation is injective;
that countably many of these subsets cover the domain of
invertibility; and that its complement -- the domain of
non-invertibility -- is measurable and so may be broken off
and handled separately.}

\vskip.10in


\noindent A.M.S. Subj. Classification 26B12, 26B15, 28A75.

\newpage

The change of variable formula in multi-dimensional calculus,
$\int_{F(E)}\varphi\,dy = \int_E\varphi(F)|J_F|\,dz$, is usually
established on the assumption that $F$ is a one-to-one selfmap of
$n$-space continuously differentiable with non-vanishing Jacobian
$J_F$
in an open domain containing $E$, $F(E)$ is a closed bounded
domain and
$\varphi$ is real-valued integrable. The integrals make sense however
-- and thus the equality can be formulated -- under more general
conditions, the natural most general being just measurability of the
domains and (Lebesgue) integrability of the integrands. Given the
standard properties of the Lebesgue integral, it suffices to obtain
the equality for the integrand $\varphi = 1$:
\begin{equation*}
|F(E) | = \int_E\Delta F'\,dz,
\end{equation*}
a formula which will be shown valid on measurable subsets $E$
on which $F$ is injective. The absolute value of the
Jacobian of $F$  has been rewritten $\Delta F'$ (for reasons
which will appear below).

If one wishes to cumulate these
subsets (which will be shown to cover $\Delta F'>0$ countably) one will need
to
replace the value of the measure on the left with the integral of the
Banach indicatrix $N(y)$, which counts the number of elements of $E$
sent by $F$ on $y$. Also, the zero-set of $\Delta F'$ will be shown to have
measure zero, which permits establishing the equality over it seperately.

This form has been attained gradually -- see the successive editions
of
Rudin -- and appears finally in his most recent edition as well as in
Smith with, however, some superfluous restrictions and -- what is
perhaps more irritating -- an appeal in both cases to Brouwer theory.
 A
purely analytic proof can be extracted by specialization from a proof
of the ``area formula'' in Evans \emph{et al}. but as given there
contains some gaps and unnecessary detours. Thus the proof below,
accessible to anyone with a modest command of Lebesgue theory, should
offer some appeal.

This paper was submitted to the Monthly in September 1998, thus was
contemporaneous with Lax's which (although far less complete) appeared in
the Monthly just when this one was rejected; a couple of years later, Lax
published an extension (II in the References) which still does not tell the
full story.



The sequel is organized as follows: The next five paragraphs
motivate and introduce the ``scale factor'', which is the
ratio by which a linear selfmap augments (or reduces)
volume, and verifies its continuity for the operator norm.
The following three paragraphs recall the definition of the
derivative for vector-valued functions and show its Borel
measurability as a function from the domain of
differentibility to the operator normed linear selfmaps. The
remaining three paragraphs then justify the integrability of
the integrand, the break-off of its zero set, and the
desired equality.

For real-valued functions of a single real variable, the
derivative equals simultaneously two {\it a priori} distinct
characteristics of a function $F$ in the neighborhood of an
argument $x$: it is the coefficient $F^\prime(x)$ of the
closest linear approximation to
$F$ near $x$ and it is the slope of the tangent to the
graph at $x.$ These entities need to be distinguished in higher
dimension; let's start with ``slope.''

For a linear function on the line, the absolute value of its slope is
 the
ratio of the length of the image of an interval by its length: it is
 the
same for all

\noindent intervals and so may as well be defined as the length of the
image of the unit semi-open interval $[0,1).$

For a linear selfmap $L$ of $n$-space, define analogously the
``slope'' or ``scale factor'' $\Delta L$ as the ratio of the volume of
the image of a box (i.e. a product of intervals) by its volume. This
will be shown in a moment to be independent of the box and so may as
well be defined as the volume of the image of the unit semi-open box
(s.o.b.), i.e. the set of $x$'s whose co-ordinates satisfy
$0\le x_i<1.$
(In the plane the image of the unit square is the parallelogram
spanned
by the images of the co-ordinate vectors, whose cross product gives
its
area; to recognize $\Delta L$ in general as the absolute value of
the determinant of the image vectors, one should arrange to obtain
these as
non--negative multiples of an orthonormal basis. However, this
identification is not
used in what follows.)

Divide the unit s.o.b.~into disjoint equal size s.o.b.'s of side
$\frac{1}{k}:$ their volume is $\frac{1}{k^n}$ and there are $k^n$ of
them, sent into each other by translation. If $L$ is invertible it
sends them into disjoint parallelepipeds, also sent into each other by
translation, hence of equal volume, whose union has volume
$\Delta L$ -- thus each image has volume
$\left(\frac{1}{k^n}\right)\Delta L.$ (If $L$ is not invertible it
maps
into a lower dimensional space and all images have $n$-volume zero.)
Any open $U$ is a disjoint union of countably many s.o.b.'s, so its
Lebesgue $n$-volume is also increased $\Delta L$ times by $L$ and thus
finally $|L(E)|=(\Delta L)|E|$ for every set $E$.

$\Delta$ is a continuous function on the $L$, topologized by
uniform convergence on bounded sets -- i.e. on any bounded
with non-void interior. The image (by every $L$) of the closed unit box
is compact, hence its epsilon neighborhoods are a base around the
image: so a uniformly convergent sequence has image measures with lim
sup $\leq$ the image measure of the limit $L$. To obtain the dual
inequality, apply the same reasoning to the complement of the open unit
box in some larger closed box.

A vector-valued function of a real variable can be differentiated in
the classical way by taking the vector limit of the difference
quotient, but this is already impossible even for a scalar-valued
function of a vector variable. What one can do is to fix a vector $u$
and take the limit of $\frac{F(z+hu)-F(z)}{h}$ as $h\downarrow 0.$
This limit would then exist for, and be positive homogeneous in,
every positive multiple $pu$ in place of $u$.
It thus suffices to postulate these limits for
unit vectors $u$; they are called the directional derivatives (of $F$
at $z$) in the directions $u$ (the unit vectors are construed to point
in the direction from the origin to their endpoint on the unit sphere;
for a function of a single real variable this is the one-sided
derivative and for co-ordinate $u$, the one-sided partials).

Assuming the directional derivatives exist in all directions $u$ at
$z,$ one could define the ``derivative'' $F^\prime(z)$ of $F$ at $z$
as
the function which takes each $u$ to the directional derivative
$F^\prime(z)u.$ One will however require the approach of
$\frac{F(z+hu)-F(z)}{h}$ to $F^\prime(z)u$ to be uniform in the unit
vectors $u$; and if $F'(z)$ is bounded this will result in
$\frac{F(x)-F(z)}{|x-z|}$ bounded for $x$ sufficiently close to $z.$
This last may be shown to entail that $F^\prime(z)pu:=pF^\prime(z)u,\
p\ge 0$, acts linearly on vectors $pu$ for a.e.~$z$: accordingly, we
now strengthen the definition by making $F^\prime(z)$ a linear
operator
whose value at every $v$ is the limit as $h \rightarrow 0$,
uniform for bounded vector
arguments $v$, of $\frac{F(z+hv)-F(z)}{h}$.

The derivative is thus a function from the domain, assumed measurable
(follows from $F$ continuous), of
differentiability of $F$ to the space of linear selfmaps of
$n$-space topologized by uniform convergence on the unit ball\footnote{This
topology is defined by the norm
$\|L\|$, the radius of the smallest ball centered at $0$ which contains
the image of the centered ball of radius one;
by linearity this $\|L\|$ is also the Lipschitz constant.}
-- a base at $0$ is obtained by taking for every centered open ball the
subset of those which send the unit ball into it. This function will
next be shown to be (Borel) measurable: i.e., the inverse image of
every open is measurable. Let $r_m\uparrow r,$ the radius of
centered open ball
$B$ in $R^n,$ which is an increasing union of closed balls $B_m$ of
radius $r_m$ with center zero. Then $F^\prime(z)u \in B$ entails (by
compactness of the $u$'s) that all but finitely many (in $k$) of the
continuous $k\left[F\left(z+\frac{u}{k}\right)-F(z)\right]$
belong to some $B_m$ (recall) uniformly in $u;$ the converse already
holds if the belonging is uniform in some countable dense subset of the
unit ball. Hence $\{z: F'(z)u\in B\}$ is the countable union (over $m$)
of the sets where a countable intersection (over $u$) of lim sup (over
$k$) of the modulus of continuous functions of $z$ and $u$ is $\leq
r_m.$

Since $\Delta L$ has been shown continuous for linear operators
$L$, it
follows that the composite, $\Delta F^\prime(z)$, is a non-negative
measurable function on a measurable subset
$E$ of $n$-space, hence has a well-defined Lebesgue integral
$\int_E \Delta F'(z)dz.$

(It also follows that the set where $\Delta F'(z) = 0$ is measurable;
since the equality for this set can be established by a simple
calculation, we shall restrict to $\Delta F'(z) > 0$ in the sequel.)

If $F':=F'(z)$ is invertible at $z$ then for $|v|<$ some
$\delta,
|F(z+v)-F(z)-F'v| \leq \bigl[\varepsilon/\|F^{\prime-1}\|\bigr]|v|
\leq \varepsilon|F'v|$; thus
$(1-\varepsilon)|F'x-F'z| \leq |Fx-Fz| \leq
(1+\varepsilon)|F'x-F'z|$ for $|x-z| < \delta$.
Observe that for fixed
$\delta$ the set of $z$'s satisfying (the extremes of)
the first inequality is measurable:
both terms are continuous in $v$ and
measurable in $z$
so imposing it at a countable dense subset of the ball of radius
$\delta$
reduces the implication to a countable set of inequalities between
measurable functions.  Since $E$ is a countable union of these sets
for a sequence
of $\delta$'s, it suffices to establish the
theorem on each of them. Then one can replace the fixed argument $z$ of
$F'$ with any $y$ in this subset and still obtain (for this $\delta$)
the extremes of the first, hence also the second, inequality with $y$
in place of $z$. 
Restricting further to measurable subsets of diameter
$<$ $\delta$ on which the values of  $|F'|$ lie between
$(1 - \varepsilon)$
and $(1 + \varepsilon)$ of its value at $z$ for every vector in the
ball of radius $\delta$, pass from the inequality with $y$ in
place of $z$ back to $F'$ at z, to attain $(1 -
\varepsilon)^2 | F'x - F' y |
\leq |Fx - Fy| \leq (1 + \varepsilon)^2 |F'x -
F'y|$ even for $x,y \neq z$.
(Since there is a countable cover by such
sets, it suffices to restrict attention to each of them.)
By the left inequality and invertibility of $F'$ at $z$,
$F$ is one-to-one on this subset.
From $|Fx-Fy|
\leq |Gx-Gy|$ for any one-to-one $F$ and all $x,y$ in some measurable
set follows: $FG^{-1}$
is non-expansive -- Lip(1) -- hence does not increase measure;
so
$|FE| \leq |FG^{-1}GE| \leq |GE|$ for each of its subsets $E$:
thus
$(1-\varepsilon)^{2n}\Delta \, F'|E|=(1-\varepsilon)^{2n}|F'(E)| \,
\leq \, |F(E)|
\,
\leq \, ( 1 + \varepsilon )^{2n}|F'(E)|= (1 +
\varepsilon )^{2n}\Delta F' |E|$.
Since $F$ is one-to-one, it preserves disjointness of subsets
and since differentiable, measurability of subsets;
hence $|F( \, \, )|$
is additive on measurable decompositions.
By passing to the limit in the definition of the integral one is led
at last to
\begin{equation*}
(1-\varepsilon)^{2n} \int_E\Delta F'\,dz \leq | F(E) | \leq
(1+\varepsilon)^{2n} \int_E\Delta F'\,dz
\end{equation*}
with arbitrary $\varepsilon$.

\vskip.10in

\bibliographystyle{amsplain}

\centerline{References}

\vskip.10in

\let\bibitem\relax


\bibitem{}
T.~M. Apostol, \emph{Mathematical Analysis}, Addison-Wesley
(1957).

\bibitem{} P.~Billingsley, \textit{Probability and Measure}, John
Wiley \& Sons (1979).



\bibitem{} L.~C. Evans \& R.~F. Gariepy, \emph{Measure Theory and Fine
Properties of Functions}, CRC Press (1992).

\bibitem{} P.~D. Lax, \emph{Change of variables in multiple integrals},
Amer. Math. Monthly, {\bf 106} (1999), 497--501; II, Ibid {\bf 108} (2001)
115--119.

\bibitem{} I.~P. Natanson, \emph{Theory of Functions of a Real
Variable}, F.~Ungar (1955).

\bibitem{} W.~F. Pfeffer, \textit{The Riemann Approach to
Integration}, Cambridge U. Press (1993).

\bibitem{} W.~Rudin, \textit{Real and Complex Analysis}, 3rd edition,
McGraw-Hill (1987).

\bibitem{} K.~T. Smith, \emph{Primer of Modern Analysis}, Springer
(1983).

\bibitem{} T.~Traynor, \emph{Change of Variable for Hausdorff Measure
(from the beginning)}, Rendiconti dell'Istit di Mat., Trieste {\bf 26}
(1994), suppl.,328--347/

\vskip.20in

\centerline {Centre de recherches math\'ematiques}
\centerline{Universit\'e de Montr\'eal}
\centerline {C.P. 6128, succursale centre-ville}
\centerline {Montr\'eal, Qc H3C 3J7}
\centerline {email:fleischi@crm.umontreal.ca}

\vskip.20in

\end{document}